
\def\Z_2{{\bf Z}_2}
\def\Z{{\bf Z}}
\def\C{{\bf C}}
\def\>{\rangle}
\def\U{$U_q[osp(1/2)]$}
\def\\U{$U_q[osp(1/2)]\;$}

 \baselineskip=16pt

\noindent
{\bf NEW SOLUTIONS OF THE YANG-BAXTER EQUATION BASED ON
ROOT OF 1 REPRESENTATIONS OF THE PARA-BOSE SUPERALGEBRA
U$_q$[osp(1/2)] }

\vskip 32pt
\noindent
T. D. Palev{*} and N. I. Stoilova{*}

\noindent
International Centre for Theoretical Physics, 34100 Trieste, Italy
\vskip 12pt

\footnote{*}{Permanent address: Institute for Nuclear Research and
Nuclear Energy, 1784 Sofia, Bulgaria; E-mail:
palev@bgearn.bitnet, stoilova@bgearn.bitnet}

\vskip 48pt

{\bf Abstract.} New solutions of the quantum Yang-Baxter equation,
depending in general on three arbitrary parameters, are written
down. They are based on the root of unity representations of the
quantum orthosymplectic superalgebra \\U, which were found recently.
Representations of the braid group $B_N$ are defined within any
$N^{th}$ tensorial power of root of 1 \\U modules.

\vfill\eject

\noindent
{\bf 1. Introduction}
\bigskip

\noindent

In the present paper we write down new solutions of the quantum
Yang-Baxter equation (QYBE), associated with root of unity
representations of the quantum orthosymplectic superalgebra
\U, which we have recently constructed [1]. All such
representations are with a highest and with a lowest weight. For
$q$ being a $4k$ root of 1 with $k=3,5,7,\ldots,$ there exists a
continuous class of $k$-dimensional representations. The
solutions of the QYBE we find depend in general on three
continuous parameters.

The general interest for studying solutions of the quantum
Yang-Baxter equation is inspired from the various applications
of the latter in conformal field theory [2, 3], quantum
integrable models [4, 5] and knot theory [6-8].  Our motivation
for the present investigation is of somewhat different nature.
It originates from the close connection between the
representations of the orthosymplectic superalgebras and the
quantum statistics [9, 10], more precisely - the parastatistics
[11].

It is perhaps worth commenting the last point in some more
details. To this end  consider as an example the Hopf algebra
$U_q[osp(1/2n)]$, the quantized universal enveloping algebra of
the orthosymplectic Lie superalgebra $osp(1/2n)$. The
quantization of the latter  in terms of its Chevalley generators
is well known [12-17]. Recently there has been given an
alternative definition of $U_q[osp(1/2n)]$ [18-21] in terms of
preoscillator generators $a_i^\pm,\;K_i=q^{H_i}$,
$i=1,\ldots,n$. The relation to the quantum statistics stems from
the observation that the operators $a_i^\pm ,\;i=1,\ldots,n$ can
be identified with deformed para-Bose operators. Moreover it
turns out that the oscillator (or Weyl) superalgebra $W_q(n)$
generated by $n$ pairs of deformed Bose operators [22-25]  is a
factor algebra of $U_q[osp(1/2n)]$ [26, 27] and (depending on the
precise definition of the preoscillator generators) a morphism
of $U_q[osp(1/2n)]$ onto $W_q(n)$ is given essentially by a
replacement of the deformed para-Bose operators with deformed
Bose operators. Therefore, despite of the fact that the
oscillator algebra $W_q(n)$ is not a Hopf algebra, one can define
an $R$-matrix associated with $W_q(n)$ simply by considering the
Fock representation of $W_q(n)$ also as a representation of
$U_q[osp(1/2n)]$. To this end one has to express the
$U_q[osp(1/2n)]$ universal $R$-matrix in terms of preoscillator
generators and subsequently replace them with deformed Bose
operators. The related matrices $R_{12},\;R_{13},\;R_{23}$,
which are functions on $n$ pairs of deformed Bose operators and
the corresponding number operators, provide a "bosonic" solution
of the QYBE. Certainly one can try to carry out the above
programme in a more general framework, considering other
representations of the preoscillator generators. This would
correspond to finding representations of the deformed para-Bose
operators.  The problem however is not simple; it has not been
solved so far even in the nondeformed case.

The present paper is a small step towards the realization of the
above programme. Here we deal with the superalgebra
$U_q[osp(1/2)]$. Nevertheless already in this simple case one
arrives to interesting conclusions.  It turns out, for instance,
that apart from the representations corresponding to both
deformed and nondeformed parabosons (and in particular - bosons)
one finds a (root of 1) representation with $a^\pm$ being usual
fermions [28], i.e., the fermions are deformed parabosons. Thus,
the bosons and the fermions appear as different irreps of one and
the same quantized superalgebra, namely \U. As an example we
write down the corresponding 4-dimensional (nondiagonal)
$R$-matrix, which leads to a "fermionic" solution of the
QYBE.

The new solutions of the QYBE will be based on
the representations of the preoscillator generators $a^\pm, \;
K=q^H$ in (deformed para-Bose) Fock spaces. We pay special
attention to the case when the deformation parameter $q$ is a
root of unity, which leads to finite-dimensional Fock spaces.

For $n>1$ the preoscillator generators of $U_q[osp(1/2n)]$ are
very different from its Chevalley generators. In case $n=1$
however  the creation and the annihilation (deformed para-Bose)
operators $a^+, \; a^-$ can be identified with the positive
and the negative root vectors $e$ and $f$ of \U, respectively.
Therefore the results to follow could have been given entirely in
terms of the canonical for \U $\;$  terminology and notation. We
prefer however to keep close to the notation of the preoscillator
generators, speaking about creation and annihilation operators
instead of Chevalley generators, (deformed) Fock spaces instead
of Verma modules, etc.  In order to underline that  \\U  is
(essentially) generated by deformed para-Bose operators we call
it a (deformed) para-Bose superalgebra.

The paper is organized as follows. In Sec. 2 we recall the
definition of the deformed para-Bose superalgebra and its Fock
irreps  for $q$ being root of 1. The form of the transformation
relations is new and more compact, comparing to those given in
[1]. In Sec. 3 new solutions of the QYBE are constructed. The
situation here is rather peculiar. We first prove that at $q$
being root of 1 \\U is in general not almost cocommutative.
Nevertheless the expression of the (generic) universal $R$-matrix
turns to be well defined within all of our representation spaces,
which leads to solutions of the QYBE. In addition the $R$-matrix
allows us to define representations of the braid group $B_N$ in
the $N^{th}$ tensorial power of any of the \\U Fock modules.

Throughout we use the following abbreviations and notation:
${\bf C}$ -- all complex numbers;
${\bf Z}$ -- all integers;
${\bf Z}_+$ -- all nonnegative integers;
${\bf Z}_2 = \{\bar{0},\bar{1}\}$;
$[A,B]=AB-BA, \hskip 4pt \{A,B\}=AB+BA$;
$U_q\equiv U_q[osp(1/2)]$.


\vskip 32pt
\noindent
{\bf 2. The para-Bose algebra U$_q$[osp(1/2)] and its Fock irreps }

\bigskip
\noindent
Here we summarize the results of [1].  The form of the
expressions (2.5), (2.7)-(2.9), describing the transformations of
the Fock spaces, is however new. It is more compact than the
corresponding relations in [1].

The superalgebra $U_q=U_q[osp(1/2)],\;q\in{\bf C} \backslash
\{0,\pm 1\}$ has three generators $a^+,\; a^-,\; H$, satisfying
the defining relations:
$$
[H,a^\pm]=\pm 2a^\pm,\quad \{a^+,a^-\}={{q^H-q^{-H}}\over{q-q^{-1}}}.
\eqno (2.1)
$$
$H$ is an even generator,  $a^\pm$ are odd.
As $q\rightarrow 1$ $H=\{a^+,a^-\}$ and eqs. (2.1) reduce to
the defining relations of the nondeformed para-Bose operators
[11] ($\xi, \eta, \epsilon =\pm \; {\rm or} \; \pm 1 $):
$$
[\{a^\xi,a^\eta\},a^\epsilon]=(\epsilon - \eta)a^\xi +(\epsilon -
\xi)a^\eta .\eqno(2.2)
$$
The Hopf algebra structure on $U_q$ can be defined in different
ways [17]. For the comultiplication we set
$$
\Delta (H)=H\otimes 1+1\otimes H,\quad
\Delta (a^+)=a^+\otimes 1+q^{-H}\otimes a^+,\quad
\Delta (a^-)=a^-\otimes q^H+1\otimes a^-. \eqno(2.3)
$$

Passing to the representations of $U_q$ we note that
the finite-dimensional irreps of $U_q[osp(1/2)]$ at generic $q$
were constructed in [29, 30].  Some root of unity highest weight
irreps were also obtained in [30]; both highest weight and cyclic
representations were studied in [31-34].

A (deformed) Fock space $F(p)$ is defined in the usual for the
parastatistics way [11]: for any complex $p$ (which is an
analogue of the  order of the parastatistics) one postulates the
existence of a vacuum vector $\vert 0\rangle \in F(p)$ so that
$a^-\vert 0\rangle=0$ and $H\vert 0\rangle=p\vert 0\rangle$.
{}From now on we shall denote by $a_p^\pm $ and $H_p$ the
representatives of $a^\pm $ and $H$ in $F(p)$. The latter
is an infinite-dimensional linear space with a basis
$\vert n\rangle=(a_p^+)^n\vert 0\rangle,\; n\in {\bf Z}_+$.

Setting
$$
\{n;x\}_q={{q^{n+x}-(-1)^n q^{-n-x}}\over{q-(-1)^n q^{-1}}} \eqno(2.4)
$$
\noindent
one can write the transformation of the basis as follows:
$$
H_p\vert n\rangle=(2n+p)\vert n\rangle ,\quad
a_p^-\vert n\rangle=\{n;0\}_q\{n-1;p\}_q\vert n-1\rangle ,\quad
 a_p^+\vert n\rangle=\vert n+1\rangle .
\eqno(2.5)
$$
At generic $q$ the space $F(p)$  is infinite-dimensional.
It is a simple (=irreducible) $U_q$ module if  $p$ is not
a negative even number [28] (which we always assume). The space
$F(p=1)$ is the Fock space of deformed Bose operators [22-25].
Within  $F(1)$ the preoscillator operators satisfy the relations
$$
a_1^-a_1^+-q^{\pm 2}a_1^+a_1^-=q^{\mp 2N}, \quad {\rm where} \;
N={1\over 2}(H_1-1)\;{\rm is\; the\; number\; operator.} \eqno(2.6)
$$
In the root of unity cases $F(p)$ is indecomposible if and only if
$q=e^{i{\pi \over 2}{m\over k}}$ for every $m,k \in {\bf Z}$ such
that  $q \notin \{\pm1 , \pm i\}$.  The factor-space of $F(p)$
with respect to the maximal invariant subspace is an
irreducible module, containing the vacuum vector $\vert 0\rangle$.

The algebras $U_q$ corresponding to all possible values of $m$
and $k$ contain several isomorphic copies. Without loss of
generality we restrict $m$ and $k$ to values, which we call
admissible, namely: (1) $ k=2,3,\ldots$; (2) $ m\in
\{1,2,\ldots,k-1\}$; (3) $m$  and $k$ are relatively co-prime.
{}From now on we consider $q=e^{i{\pi \over 2}{m\over k}}$  to be
only an admissible root of 1.

The irreducible $U_q$ modules with $q$ being root of 1
are finite-dimensional. Denote by $W^L(p) \subset F(p)$ an
$L+1$-dimensional representation space with a basis
$\vert 0\rangle,\;\vert 1\rangle,\ldots, \vert L\rangle$.
Its transformations under the action of the $U_q$ generators read:
$$
H_p\vert n\rangle=(2n+p)\vert n\rangle ,\;
a_p^-\vert n\rangle=\{n;0\}_q\{n-1;p\}_q\vert n-1\rangle ,\;
a_p^+\vert L\rangle=0,\; a_p^+\vert n\rangle=\vert n+1\rangle ,\;
n<L. \eqno(2.7)
$$
We distinguish two classes of algebras, each one containing
three groups of representations,
\noindent
$$
\eqalignno{
& Class\; I\;(k-m={\rm odd}): \;(I.a)\;L=2k-1 \;{\rm if}\;
  p\neq {\rm integer}; \;  (I.b)\;L=p(k-1)(mod\; 2k) \;
  {\rm if}\; p={\rm  integer};\cr
& \hskip 39mm   (I.c)\; L=2k-1 & (2.8) \cr
& Class\; II \;(k,m={\rm odd}):\;\; (II.a)\; L=k-1
  \;{\rm if} \; p\not={\rm even};  \;
  (II.b)\;L=(k-p)(mod\;k) \;{\rm if}\; p={\rm even}\cr
& \hskip 39mm (II.c)\;L=k-1. & (2.9)\cr
}
$$
\noindent
The cases $(I.a),\;(I.b),\;(II.a)$ and $(II.b)$ correspond to
irreducible representations, whereas in $(I.c)$ (resp.  $(II.c)$)
the representation is indecomposible if $p={\rm integer}$ (resp.
if $p={\rm even}$).  The $2k$-dimensional modules corresponding
to $(I.c)$  were described in [32], where in particular it was
shown how one can modify those of them  corresponding to $k$=odd
and $m$=even so that they carry cyclic representations. One has
to keep in mind however that at certain values of $p$ these
modules are no more irreducible, they are indecomposible.  In
fact each simple module $W^L(p)$ from $(I.b)$ with
$L=p(k-1)(mod\; 2k)$ is a factor-module of $W^{2k-1}(p)$ from
$(I.c)$ with respect to its maximal invariant subspace. To our
best knowledge the representations from the classes $(I.b)$ and
$II$ were not described in the literature so far.

One can always assume $0<Re(p)\leq 4k$, since the representations
with $p$ outside that interval are equivalent to representations
with $p$ obeying the above inequality; if $k$ is odd and $m$ is
even  one can further set $0<Re(p)\leq 2k$ if $m=2(mod\;4)$ and
$0<Re(p)\leq k$ if $m=4(mod\;4)$.


\vskip 32pt
\noindent
{\bf 3. R-matrices and new solutions of the QYBE }

\bigskip
\noindent
One way for constructing $R$-matrices and hence solutions of the
QYBE is based on the use of the universal $R$-matrix of a
quasitriangular Hopf algebra $U$ together with the
representations of $U$.

The universal $R$-matrix for $U_q$ was written down in [29, 30].
Here we use the expression as given in [17], which read in our
notation:
$$
R=\sum_{n\geq 0}(-1)^{n(n+1)\over 2}{(q-\bar q)^n\over
{{(n)_{-\bar q^2}!}}}[(a^+)^n\otimes (a^-)^n]q^{{1\over 2}
H\otimes H} \eqno(3.1)
$$
where $\bar q=q^{-1},\;\; (n)_a=(1-a^n)/(1-a),\;\;
(n)_a!=(1)_a(2)_a\ldots(n)_a. $

If $\rho_1$ and $\rho_2$ are two representations of $U_q$
defined in $V_1$ and $V_2$, then the related
$R$-matrix is $R(\rho_1,\rho_2)=(\rho_1\otimes \rho_2)R\in
End(V_1\otimes V_2)$.

In the root of 1 cases however the above construction generally
fails , because for certain admissible $q$ $U_q$ is no more
almost cocommutative.  The proof is essentially the same as the
one given  by Arnaudon for $U_q[sl(2)]$ [35]. It is based on the
observation that  $U_q$ contains a larger center generated from
its Casimir operator and the additional central elements
$\hat{x}^\pm=(a^\pm)^{2k}$ and $\hat{z}=(K)^{2k}$ [30, 31]. If
$\rho$ is an irrep of $U_q$ in $V$, then
$\rho(\hat{x}^\pm)=\rho(x^\pm) {\bf 1}_V$,
$\rho(\hat{z})=\rho(z) {\bf 1}_V$,
where ${\bf 1}_V$ is the unit operator in $V$ and
$\rho(x^\pm),\;\rho(z)\in \C$.

We proceed to show that the universal $R$-matrix does not exist
for a subclass of I, corresponding to all algebras with $k$=odd
and $m$=even. To this end we use the following general identity:
if $q$ is an $l^{th}$ primitive root of 1 and $AB+q^2BA=0$, then
for any $N\leq l$
$$
(A+B)^N=\sum_{n=0}^N q^{-n(N-n)}{N\brace n}_q A^n B^{N-n},
\quad {N\brace n}_q={\{N\}_q!\over \{n\}_q!\{N-n\}_q!},
\quad  \{n\}_q=q^n-(-1)^nq^{-n}.
\eqno(3.2)
$$
Applying (3.2) for $N=2k$, $A=1\otimes a^-$ and $B=a^-\otimes K$,
we obtain for all class I algebras:
$$
\Delta(\hat{x}^-)=1\otimes \hat{x}^- + \hat{x}^-\otimes \hat{z},
\quad \Delta^{op}(\hat{x}^-)=\hat{x}^-\otimes 1 +
\hat{z}\otimes \hat{x}^- .
\eqno(3.3)
$$
In (3.3) $\Delta^{op}=\sigma \Delta$ is the opposite
comultiplication; $\sigma$ is a superpermutation,
$\sigma(a\otimes b)=(-1)^{deg(a)deg(b)}b\otimes a$.  Assume now
that $U_q$, is almost cocommutative, namely that there exists an
invertible element $R$ from (the completion of) $U_q \otimes
U_q$, so that $R\Delta(a)=\Delta^{op}(a)R$ for any $a\in U_q$.
On the tensor product of two irreps $\rho_1$ and $\rho_2$ in
$V_1$ and $V_2$ one would have for $a=\hat{x}^-$: $$
(\rho_1\otimes \rho_2)(R)
(\rho_1\otimes \rho_2)(\Delta(\hat{x}^-))=
(\rho_1\otimes \rho_2)(\Delta^{op}(\hat{x}^-))
(\rho_1\otimes \rho_2)(R).\eqno(3.4)
$$
Acting with both sides of (3.4) on an arbitrary vector
$\vert X\rangle \in V_1\otimes V_2$, one gets:
$$
\{\rho_1(x^-)+\rho_1(z)\rho_2(x^-)\}\vert Y\rangle=
\{\rho_2(x^-)+\rho_2(z)\rho_1(x^-)\}\vert Y\rangle ,
$$
where $\vert Y\rangle = (\rho_1\otimes \rho_2)(R) \vert X\rangle$.
Therefore
$$
\rho_2(x^-)+\rho_2(z)\rho_1(x^-)=
\rho_1(x^-)+\rho_1(z)\rho_2(x^-). \eqno(3.5)
$$
The central elements $\hat{x}^\pm$ and $\hat{z}$ can take arbitrary
values on the cyclic irreps of the algebras with $k$=odd and
$m$=even [34], i.e., in this case  $\rho_1(x^-)$, $\rho_2(x^-)$,
$\rho_1(z)$ and $\rho_2(z)$ are arbitrary numbers, which is in
contradiction to (3.5). Therefore the universal $R$-matrix cannot
exists for these algebras.  Note however that (3.5) is not in
contradiction with the representations (2.8), since for any of
them $\rho(x^\pm)=0$. Therefore, following Rosso [36],
one can try to produce an almost universal $R$-matrix on the
quotient $\tilde{U}_q=U_q[osp(1/2)]/(\hat{x}^\pm=0)$.

Eq. (3.5) holds for the subclass of the class I algebras,
corresponding to $k$=even and $m$=odd. The known
irreps for this subclass are only those listed
in (2.8).  The latter do not contradict (3.5), since
$\hat{x}^\pm$ and $\hat{z}$ act as zero operators within each
class I $U_q$-module.  Therefore the question whether the
universal $R$ matrix exists for the algebras with $k$=even and
$m$=odd is an open one. The same holds for all algebras from the
class II. Within each $U_q$ module corresponding to (2.9)
$\hat{x}^\pm$ and $\hat{z}$ are zero operators. Our attempts to
extend these modules to carry cyclic representations were not
successful.  Moreover, the eqs. (3.3) and hence (3.5) are no more
true.

We see that the question about the existence of an universal $R$
matrix for the algebras $U_q$ with $q$ being root of 1 cannot be
answered uniquely at present.  Our claim is that the $R$-matrix
(3.1), considered as an element from $\tilde{U}_q\otimes
\tilde{U}_q $ is almost universal, namely it is "universal" for
all Fock representations (2.7)-(2.9): if $\rho^{L_1}(p_1)$ and
$\rho^{L_2}(p_2)$ are  any two such representations, then the
operator

$$
R^{L_1, L_2}(p_1,p_2)=(\rho^{L_1}(p_1) \otimes \rho^{L_2}(p_2))(R):\;
W^{L_1}(p_1)\otimes W^{L_2}(p_2) \rightarrow
W^{L_1}(p_1)\otimes W^{L_2}(p_2)
\eqno(3.6)
$$
satisfies the analogue of (3.4),
$$
R^{L_1, L_2}(p_1,p_2)(\rho^{L_1}(p_1) \otimes \rho^{L_2}(p_2))
(\Delta(a))=(\rho^{L_1}(p_1) \otimes \rho^{L_2}(p_2))
(\Delta^{op}(a))R^{L_1, L_2}(p_1,p_2). \eqno(3.7)
$$
The explicit action of $R^{L_1, L_2}(p_1,p_2)$ on the basis
$\vert l_1\rangle \otimes \vert l_2\rangle $ of
$W^{L_1}(p_1)\otimes W^{L_2}(p_2)$ yields:
$$
\eqalignno
{
& R^{L_1, L_2}(p_1,p_2)(\vert l_1\rangle \otimes
\vert l_2\rangle ) = q^{{1\over
2}(2l_1+p_1)(2l_2+p_2)}\sum_{n=0}^{min(L_1-l_1, l_2)}
(-1)^{{n\over 2}(n+2l_1+1)}{(q-\bar q)^n\over
{{(n)_{-\bar q^2}!}}} & \cr
& \hskip 30mm \times \prod_{i=0}^{n-1}\{l_2-i;0\}_q\{l_2-1-i;p_2\}_q
\vert l_1+n\rangle\otimes \vert l_2-n\rangle. & (3.8)\cr
}
$$
The proof of (3.7) is by a direct computation within each $U_q$
module $W^{L_1}(p_1)\otimes W^{L_2}(p_2)$, i.e., using the
transformation relations (3.8).

The linear operators
$$
R_{12}^{L_1, L_2}(p_1,p_2),
R_{13}^{L_1, L_3}(p_1,p_3),
R_{23}^{L_2, L_3}(p_2,p_3)\;{\rm in} \;
W^{L_1,L_2,L_3}(p_1,p_2, p_3)\equiv
W^{L_1}(p_1)\otimes W^{L_2}(p_2)\otimes W^{L_3}(p_3)\eqno(3.9)
$$
which satisfy the QYBE
$$
  R_{12}^{L_1, L_2}(p_1,p_2)
  R_{13}^{L_1, L_3}(p_1,p_3)
  R_{23}^{L_2, L_3}(p_2,p_3)
= R_{23}^{L_2, L_3}(p_2,p_3)
  R_{13}^{L_1, L_3}(p_1,p_3)
  R_{12}^{L_1, L_2}(p_1,p_2).  \eqno(3.10)
$$
are defined on the basis as follows:
$$
\eqalignno
{
& R_{12}^{L_1, L_2}(p_1,p_2)(\vert l_1\rangle
  \otimes \vert l_2\rangle \otimes \vert l_3\rangle )
  = q^{{1\over 2}(2l_1+p_1)(2l_2+p_2)}
  \sum_{n=0}^{min(L_1-l_1, l_2)}
  (-1)^{{n\over 2}(n+2l_1+1)} & \cr
&&\cr
& \hskip 10mm \times {(q-\bar q)^n\over{{(n)_{-\bar q^2}!}}}
  \prod_{i=0}^{n-1}\{l_2-i;0\}_q\{l_2-1-i;p_2\}_q
  \vert l_1+n\rangle\otimes \vert l_2-n\rangle
  \otimes \vert l_3\rangle. & (3.11)\cr
}
$$

$$
\eqalignno
{
& R_{13}^{L_1, L_3}(p_1,p_3)(\vert l_1\rangle
  \otimes \vert l_2\rangle \otimes \vert l_3\rangle )
  = q^{{1\over 2}(2l_1+p_1)(2l_3+p_3)}
  \sum_{n=0}^{min(L_1-l_1, l_3)}
  (-1)^{{n\over 2}(n+2l_1+2l_2+1)} & \cr
&&\cr
& \hskip 10mm \times {(q-\bar q)^n\over{{(n)_{-\bar q^2}!}}}
  \prod_{i=0}^{n-1}\{l_3-i;0\}_q\{l_3-i-1;p_3\}
  \vert l_1+n\rangle\otimes \vert l_2\rangle
  \otimes \vert l_3-n\rangle. & (3.12)\cr
}
$$

$$
\eqalignno
{
& R_{23}^{L_2, L_3}(p_2,p_3)(\vert l_1\rangle
  \otimes \vert l_2\rangle \otimes \vert l_3\rangle )
  = q^{{1\over 2}(2l_2+p_2)(2l_3+p_3)}
  \sum_{n=0}^{min(L_2-l_2, l_3)}
  (-1)^{{n\over 2}(n+2l_2+1)} & \cr
&&\cr
& \hskip 10mm \times {(q-\bar q)^n\over{{(n)_{-\bar q^2}!}}}
  \prod_{i=0}^{n-1}\{l_3-i;0\}_q\{l_3-i-1;p_3\}
  \vert l_1\rangle\otimes \vert l_2+n\rangle
  \otimes \vert l_3-n\rangle. & (3.13)\cr
}
$$
The operators (3.9) can be expressed in terms of the $R$-matrix
(3.8). To this end introduce a superpermutation linear operator
$
P_{23}:
(\vert n_1\rangle \otimes \vert n_2\rangle)\otimes \vert n_3\rangle)=
(-1)^{n_2n_3}\vert n_1\rangle \otimes \vert n_3\rangle
\otimes \vert n_2\rangle).
$
Then
$$
\eqalign{
& R_{12}^{L_1, L_2}(p_1,p_2)=R^{L_1, L_2}(p_1,p_2)\otimes 1,\cr
& R_{13}^{L_1, L_3}(p_1,p_3)=
  P_{23}(R^{L_1, L_3}(p_1,p_3)\otimes 1)P_{23}, \cr
& R_{23}^{L_2, L_3}(p_2,p_3)=
  1\otimes R^{L_2, L_3}(p_2,p_3). \cr
}\eqno(3.14)
$$

Depending on the choice of the representations (2.8) and (2.9)
one obtains $R$-matrices of different dimensions, which may be
parameter independent or can depend on one or two free parameters.

If $\rho^{L_1}(p_1),\;\rho^{L_2}(p_2)\in (I.c)$, then
$R^{L_1, L_2}(p_1,p_2)$ depends on two arbitrary complex
parameters $p_1$ and $p_2$, $dim(R^{L_1, L_2}(p_1,p_2))=4k^2$. These
$R$-matrices were obtained in [32]. The expression (3.8) is
somewhat more compact.

If $\rho^{L_1}(p_1),\;\rho^{L_2}(p_2)\in (II.c)$
$R^{L_1, L_2}(p_1,p_2)$ depends also on the arbitrary complex
parameters $p_1$ and $p_2$, but $dim(R^{L_1, L_2}(p_1,p_2))=k^2$.
This is a new class of $R$-matrices, leading through (3.14) to
new solutions of the QYBE, defined in a $k^3$-dimensional space
$W^{L_1,L_2,L_3}(p_1,p_2, p_3)$ with $k=3,5,7,\ldots$ and
depending on three arbitrary parameters.

In all other cases the $R$-matrices depend on less then two free
parameters, which is due to the case that for certain values of
$p_1,\; p_2$ and $p_3$ $W^{L_1,L_2,L_3}(p_1,p_2, p_3)$ contains
invariant subspaces.  Those, corresponding to
$\rho^{L_1}(p_1),\;\rho^{L_2}(p_2)\in (I.b)$ or $(II.b)$ lead to
constant $R$-matrices and hence to constant solutions of the QYBE.
Here are some examples.

\noindent
{\it Example 1.} The representation  $(I.b)$ with $k=2, \;
(m=1)$ and $p=1$ gives $L=1$. From (2.7) one concludes that $a^\pm$
are Fermi operators. In the basis
$\{\vert 0\rangle\otimes \vert 0\rangle,\;
\vert 0\rangle\otimes \vert 1\rangle,\;
\vert 1\rangle\otimes \vert 0\rangle,\;
\vert 1\rangle\otimes\vert 1\rangle \} $
the "fermionic" $R$-matrix reads:
$$
R^{L_1=1,L_2=1}(p_1=1,p_2=1)=\left(\matrix{
e^{i\pi\over 8} & 0 & 0 & 0 \cr
0 & e^{3i\pi\over 8} & 0 & 0 \cr
0 & e^{i\pi\over 8}-e^{5i\pi \over 8} & e^{3i\pi\over 8} & 0 \cr
0 & 0 & 0 & -e^{i\pi\over 8} \cr
}\right).\eqno(3.15)
$$
It contains no free parameters.

\noindent
{\it Example 2.} We consider the class II algebra $U_q$ with the
smallest possible value of $k$, namely $k=3$ (and hence
$m=1$), i.e., $q=e^{i\pi/6}$. There is a tree of $R$-matrices,
related to the different possible branches of the representations
(II.a,b,c).  One such branch is, for instance, $R^{2,2}(p_1,p_2)
\rightarrow R^{2,1}(p_1,2)\rightarrow R^{1,1}(2,2)$. The root
$R$-matrix $R^{L_1=2,L_2=2}(p_1,p_2)$ is 9-dimensional and
depends on two arbitrary parameters $p_1$ and $p_2$. In a matrix
form ( ordering the basis lexically, $\vert i\rangle\otimes \vert
j\rangle <
\vert k\rangle\otimes \vert l\rangle$ if $i<k$ or
if $i=k$ and $j<l$) one
obtains from (3.8):

$$
R^{2,2}(p_1,p_2)=
\left(\matrix{
A_{00,00} & 0 & 0 & 0 & 0 & 0 & 0 & 0 & 0 \cr
0 & A_{01,01} & 0 & 0 & 0 & 0 & 0 & 0 & 0 \cr
0 & 0 & A_{02,02} & 0 & 0 & 0 & 0 & 0 & 0 \cr
0 & A_{10,01} & 0 & A_{10,10} & 0 & 0 & 0 & 0 & 0 \cr
0 & 0 & A_{11,02} & 0 & A_{11,11} & 0 & 0 & 0 & 0 \cr
0 & 0 & 0 & 0 & 0 & A_{12,12} & 0 & 0 & 0 \cr
0 & 0 & A_{20,02} & 0 & A_{20,11} & 0 & A_{20,20} & 0 & 0 \cr
0 & 0 & 0 & 0 & 0 & A_{21,12} & 0 & A_{21,21} & 0 \cr
0 & 0 & 0 & 0 & 0 & 0 & 0 & 0  & A_{22,22} \cr
}\right).\eqno(3.16)
$$
with
$$
\eqalign{
& A_{00,00}=e^{{i\pi\over 12}p_1p_2}\;\;
  A_{01,01}=e^{{i\pi\over 12}p_1(p_2+2)}\;\;
  A_{02,02}=e^{{i\pi\over 12}p_1(p_2+4)}\;\;
  A_{10,10}=e^{{i\pi\over 12}(p_1+2)p_2}\;\;
  A_{11,11}=e^{{i\pi\over 12}(p_1+2)(p_2+2)}\cr
&\cr
& A_{12,12}=e^{{i\pi\over 12}(p_1+2)(p_2+4)}\;\;
  A_{20,20}=e^{{i\pi\over 12}(p_1+4)p_2}\;\;
  A_{21,21}=e^{{i\pi\over 12}(p_1+4)(p_2+2)}\;\;
  A_{22,22}=e^{{i\pi\over 12}(p_1+4)(p_2+4)} \cr
&\cr
& A_{10,01}=-2ie^{{i\pi\over 12}p_1(p_2+2)}sin({\pi\over 6}p_2)\;
  A_{11,02}=-2ie^{{i\pi\over 12}p_1(p_2+4)}cos({\pi\over 6}(p_2+1))\cr
&\cr
& A_{20,11}=2ie^{{i\pi\over 12}(p_1+2)(p_2+2)}sin({\pi\over 6}p_2)
A_{20,02}=-ie^{{i\pi\over 12}(p_1(p_2+4)+2)}
  (2sin({\pi\over 6}(2p_2+1))-1)\cr
&\cr
& A_{21,12}=2ie^{{i\pi\over 12}(p_1+2)(p_2+4)}cos({\pi\over 6}(p_2+1)).\cr
}
$$
Setting $p_2=2$  and $L_2=1$ one obtains the next matrix from the
branch, namely the
6-dimensional $R$-matrix $R^{L_1=2,L_2=1}(p_1,p_2=2)$, which depends
on the arbitrary parameter $p_1$:
$$
R^{2,1}(p_1,2)=
\left(\matrix{
e^{{i\pi\over 6}p_1} & 0  & 0 & 0 & 0 & 0 \cr
0 & e^{{i\pi\over 3}p_1} & 0 & 0 & 0 & 0 \cr
0 & -i\sqrt{3}e^{{i\pi\over 3}p_1} & e^{{i\pi\over 6}(p_1+2)}
	  & 0 & 0 & 0 \cr
0 & 0 & 0 & e^{{i\pi\over 3}(p_1+2)}  & 0 & 0 \cr
0 & 0 & 0 & i\sqrt{3}e^{{i\pi\over 3}(p_1+2)} &
	  e^{{i\pi\over 6}(p_1+4)} & 0 \cr
0 & 0 & 0 & 0 & 0 & e^{{i\pi\over 3}(p_1+4)}  \cr
}\right).\eqno(3.17)
$$
$R^{2,1}(p_1,2)$ can be obtained from the root matrix (3.16) by
crossing out its rows and columns with numbers 3, 6 and 9
and setting $p_2=2$.

The last matrix from the branch corresponds to
$p_1=p_2=2$ and $L_1=L_2=1$. It is a 4-dimensional
constant $R$-matrix, which can be obtained by crossing out the
last two rows and columns in (3.17) and setting $p_1=2$:
$$
R^{1,1}(2,2)=
\left(\matrix{
e^{i\pi\over 3} & 0  & 0 & 0  \cr
0 & e^{2i\pi\over 3} & 0 & 0 \cr
0 & -i{\sqrt{3}}e^{2i\pi\over 3} & e^{2i\pi\over 3} & 0 \cr
0 & 0 & 0 & -e^{i\pi\over 3}\cr
}\right).\eqno(3.18)
$$
One can choose certainly  other branches from the $R$-matrix
tree, obtaining in this way new $R$-matrices of smaller
dimensions, which are always submatrices of the root matrix
(3.16).

Let us mention at the end, following Zhang [38], that the
$R$-matrix can be used also in order to define representations of
the braid group $B_N$ acting in any $N^{th}$ tensorial power of
Fock spaces $W^L(p)$, namely in $W^L(p)^{\otimes N}$. To this end
set $\check{R}^L(p)=PR^{L,L}(p,p)\in End(W^L(p)\otimes W^L(p))$,
where $P$ is the superpermutation operator in $W^L(p)\otimes
W^L(p)$. It is straightforward to verify that $\check{R}^L(p)$
is an \\U intertwining operator in $W^L(p)\otimes W^L(p)$:
$$
[\check{R}^L(p),\Delta(a)]=0 \quad \forall a\in
U_q.\eqno(3.19)
$$
Hence [38] $\sigma_i \in End(W^L(p))^{\otimes N} \;\;
i=1,\ldots,N-1$, defined as
$$
\sigma_i={\bf 1}^{\otimes (i-1)}\otimes \check{R}^L(p)\otimes
{\bf 1}^{\otimes (N-i-1)} \eqno(3.20)
$$
gives a representation of $B_N$, namely the
$\sigma_1,\ldots,\sigma_{N-1}$
satisfy the defining relations for $B_N$:
$$
\sigma_i\sigma_j=\sigma_j\sigma_i \;\;\vert i-j\vert >1;\quad
\sigma_i\sigma_{i+1}\sigma_i=\sigma_{i+1}\sigma_i\sigma_{i+1}.
\eqno(3.21)
$$
Hence (the representation of the braid group) $B_N$ is a subset
of the set of all intertwining operators in $W^L(p)^{\otimes N}$.


\vskip 32pt
\noindent
{\bf 4. Concluding remarks}

\bigskip
\noindent
We have found new solutions of the quantum Yang-Baxter equations,
using essentially the representations of \\U, which we have
recently constructed.  The solutions were obtained formally from
the "generic" $R$-matrix (3.1), despite of the fact that the
latter does not exists in root of 1 cases. The more precise
statement is that at values of the deformation parameter
$q=e^{i{\pi \over 2}{m\over k}}$ with $k=$odd and $m$=even \\U is
not quasitriangular and even less - it is not almost
cocommutative. In all other admissible cases the question about
the existence of $R$ is an open one.

The results we have announced in the present paper are more
of a mathematical nature. The very fact however that $a^\pm$ are
deformed para-Bose operators (in some other terminology -
deformed supersingletons [39]) indicates already their relation
to quantum physics. In fact the representation with $p=1$
corresponds deformed to bosons [22-25]. The one-dimensional
quantum oscillator based on such operators exhibits quite unusual
properties at $q$ being root of 1. In particular it leads to
discretization of the spectrum of the position and the momentum
operators, thus putting the phase space on a lattice [40]. It
will be interesting to consider the same problem in the frame of
the more general para-Bose oscillator, considering all its
(unitarizable) root of 1 representations.

Various kinds of oscillators based on deformed parabosons were
discussed in the literature so far (see [1] for references in
this respect) without usually paying attention to the underlying
coalgebra structure. The arbitrary deformations may face however
serious problems: if the underlying deformed para-Bose algebra is
not a Hopf algebra (or at least an associative algebra with a
comultiplication, which is an algebra morphism), it is impossible
to define tensor products of representations.  The deformations
of the parabosons we consider are free of this disadvantage,
since our deformed algebra is identical with the Hopf algebra \U.
Another positive feature of the Hopf algebra deformations is the
existence of an $R$-matrix within every Fock space $W^L(p)$. The
latter allows one to define an action of the braid group $B_N$
within any $N^{th}$ tensorial power $W^L(p)^{\otimes N}$, which
commutes with \U. This is a step towards the decomposition of
$W^L(p)^{\otimes N}$ into irreducible \\U modules.

It will be interesting to generalize the present approach to the
case of several, say, $n$ modes of preoscillator operators. To
this end one has first to express the universal $U_q[osp(1/2n)]$
$R$-matrix in terms of deformed para-Bose operators and then
consider root of 1 representations of them. A good candidate for
such a representation is the one of the $q-$commuting deformed
Bose operators, introduced recently in [20, 21], which permit
only root of 1 (unitary) representations.

\vskip 24pt
\noindent
{\bf Acknowledgments}

\bigskip
\noindent
All results of the present investigation were obtained during our
three months visit at the International Centre for Theoretical
Physic in Trieste.  We are grateful to Prof. Randjbar-Daemi for
the invitation and for the kind hospitality at the High Energy
Section of ICTP.  We are thankful to Prof. D. Arnaudon for the
constructive remarks and suggestions. The work was supported by
the Grant $\Phi-416$ of the Bulgarian Foundation for Scientific
Research.

\vskip 24pt

\noindent
{\bf References}
\vskip 12pt
{\settabs \+  [11] & I. Patera, T. D. Palev, Theoretical
   interpretation of the experiments on the elastic \cr

\+ [1] & Palev T D and Stoilova N I, Unitarizable representations
	 of the deformed para-Bose \cr
\+         & superalgebra $U_q[osp(1/2)]$ at roots of 1, 1995 {\it Preprint}
	 ICTP IC/95/97 \cr

\+ [2] & Alvarez-Gaum\'e, Gomez C and Sierra G 1989 {\it Phys.
	 Lett. B} {\bf 220} 142; {\it Nucl. Phys. B} {\bf 319} 155\cr

\+ [3] & Moore G and Reshetikhin N Y 1989
		 {\it Nucl. Phys. B} {\bf 328} 557 \cr

\+ [4] & Baxter R J {\it Exactly Solved Models in Statistical
	 Mechanics} Academic Press, New York, 1982 \cr

\+ [5] & Fadeev L D Integrable models in (1+1)-dimensional
	 quantum field theory, in \cr
\+     & {\it Field Theory and Statistical Mechanics}
	 Norht-Holland, Amsterdam 1984, p.563 \cr

\+ [6] & Witten E 1989 {\it Comm. Math. Phys.} {\bf 121} 351 \cr

\+ [7] & Reshetikhin N Quantized universal enveloping algebras,
	 the Yang-Baxter equations \cr
\+         & and invariants of links: I, II {\it Preprints} LOMI
	 E-487, E-17-87 \cr
\+ [8] & Zhang R B, Gould M D and Bracken A J 1991
	 {\it Comm. Math. Phys.} {\bf 137} 13 \cr

\+ [9] & Ganchev A and Palev T D 1978 {\it Preprint} JINR P2-11941;
	  1980 {\it J. Math. Phys.} {\bf 21} 797 \cr

\+ [10] & Palev T D  1982 {\it J. Math. Phys.} {\bf 23} 1100 \cr

\+ [11] & Green H S 1994 {\it Austr.\ J.\ Phys.} {\bf 47} 109 \cr

\+ [12] & Chaichan M and Kulish P 1990 {\it Phys.\ Lett.\ B}
	  {\bf 234} 72 \cr

\+ [13] & Bracken A J, Gould M D and Zhang R B 1990
	  {\it Mod. Phys. Lett. A} {\bf 5} 331 \cr

\+ [14] & Floreanini R, Spiridonov V P and Vinet L 1990
	  {\it Phys.\ Lett.\ B} {\bf 242} 383 \cr

\+ [15] & Floreanini R, Spiridonov V P and Vinet L 1991
	  {\it Commun.\ Math.\ Phys.}{\bf 137} 149 \cr

\+ [16] & d\'{}$\!\!$ Hoker, Floreanini R and Vinet L 1991
	  {\it Journ. Math. Phys.} {\bf 32} 1247 \cr

\+ [17] & Khoroshkin S M  and Tolstoy  V N  1991 {\it Commun.\
	  Math.\ Phys.} {\bf 141} 599 \cr

\+ [18] & Palev T D 1993 {\it J. Phys. A} {\bf 26} L1111  \cr

\+ [19] & Hadjiivanov L K 1993 {\it Journ. Math. Phys.} {\bf 34}
	 5476 \cr

\+ [20] & Palev T D 1994 {\it Journ. Group Theory in Physics}
          {\bf 3} 1 \cr

\+ [21] & Palev T D and Van der Jeugt J 1995
   	      {\it J. Phys. A: Math. Gen.} {\bf 28} 2605  \cr

\+ [22] & Macfarlane A J 1989 {\it J.\ Phys.\ A~: Math.\ Gen.}
    	  {\bf 22}  4581  \cr

\+ [23] & Biedenharn L C 1989 {\it J.\ Phys.\ A~: Math.\ Gen.} {\bf 22}
    	  L873 \cr

\+ [24] & Sun C P and Fu H C 1989 {\it J.\ Phys.\ A~: Math.\ Gen.}
	     {\bf 22}  L983  \cr

\+ [25] & Hayashi T 1990 {\it Commun.\ Math.\ Phys.} {\bf 127}  129  \cr

\+ [26] & Palev T D and Stoilova N I 1993 {\it Lett. Math. Phys.} {\bf 28} 187
\cr

\+ [27] & Palev T D 1993 {\it Lett. Math. Phys.} {\bf 28} 321 \cr

\+ [28] & Celeghini E, Palev T D and Tarlini M 1990 {\it Preprint}
	 YITP/K-865 Kyoto and \cr
\+      & 1991  {\it Mod. Phys. Lett. B} {\bf 5} 187    \cr

\+ [29] & Kulish P P and Reshetikhin N Yu 1989 {\it Lett. Math. Phys}
	  {\bf 18} 143 \cr

\+ [30] & Saleur H 1990 {\it Nucl. Phys. B} {\bf 336} 363 \cr

\+ [31] & Sun Chang-Pu, Fu Hong-Chen and Ge Mo-Lin 1991 {\it
	  Lett. Math. Phys} {\bf 23} 19 \cr

\+ [32] & Ge Mo-Lin, Sun Chang-Pu and Xue Kang 1992 {\it Phys.
	  Lett. A} {\bf 163} 176 \cr

\+ [33] & Sun Chang-Pu 1993 {\it N. Cim. B} {\bf 108} 499 \cr

\+ [34] & Fu Hong-Chen and Ge Mo-Lin 1992 {\it Commun. Theor. Phys.}
	  {\bf 18} 373 \cr

\+ [35] & Arnaudon D 1995 {\it Lett. Math. Phys}
	  {\bf 33} 39 \cr

\+ [36] & Rosso M in 1992 {\it Topological and Geometrical Methods
	      in Field Theory} \cr
\+      & (Editors Mickelsson J and Pekonen O, World Scientific)
           p. 347 \cr

\+ [37] & Palev T D and Van der Jeugt J 1995
	  {\it J.\ Phys.\ A~: Math.\ Gen.} {\bf 28} 2605 \cr

\+ [38]  & Zhang R B 1992 {\it Journ. Math. Phys.} {\bf 33} 3918 \cr

\+ [39]  & Flato M, Hadjiivanov L K and Todorov I T 1993
          {\it Found. Phys.} {\bf 23} 571 \cr

\+ [40] & Bonatsos D, Daskaloyannis C, Demosthenes E and
          Faessler A 1994 {\it Phys. Lett.} {\bf B331} \cr

\end